\documentclass[preprint,12pt]{elsarticle}

\usepackage{graphicx,epsfig}

\usepackage{amssymb,amsmath}
\usepackage{graphicx}

\journal{}
\begin{document}
\begin{frontmatter}
\newtheorem{thm}{Theorem}[section]
\newtheorem{lem}{Lemma}[section]
\newdefinition{definition}{Definition}
\newdefinition{rmk}{Remark}
\newproof{pf}{Proof}
\newtheorem{ex}{Example}
\title{A Method for Unknotting Torus Knots\tnoteref{t1}}
\tnotetext[t1]{This document is a collaborative effort.}

\author[v]{Vikash. S}

\author[mp]{Madeti. P\corref{cor1}}
\ead{prabhakar@iitrpr.ac.in}
\ead[url]{http://www.iitrpr.ac.in/html/faculty/prabhakar.shtml}
\cortext[cor1]{Corresponding author}
\address[v]{Department of Mathematics, IIT Ropar, Rupnagar- 140001, India.}
\address[mp]{Room No. 207, Department of Mathematics, IIT Ropar, Rupnagar - 140001, India.}

\begin{abstract}
Unknotting numbers for torus knots and links are well known. In this paper, we present a method for determining the position of unknotting number crossing changes in a toric braid $B(p,q)$ such that the closure of the resultant braid is equivalent to the trivial knot or link. Also, we provide a simple proof for the important result $K(p,q)\sim K(q,p)$ using the results obtained from this method.
\end{abstract}

\end{frontmatter}


\section{Introduction}

The unknotting number $u(K)$ of a knot $K$ is the minimum number of crossing changes required to convert $K$ into a
trivial knot taken over all knot diagrams representing $K$. An unknotting operation involves switching the under/over crossing
strand at a crossing point of a knot diagram. One important tool to represent knots is that of the closure of braids. An $n$-string
braid is a set of $n$ smooth non-intersecting and non-self intersecting strings in $\mathbb{R}^3$, whose initial points
are $(i, 0, 0), \ i = 1,\ldots n,$ and final points are $(j, 0, 1), \ j = 1,\ldots n,$ such that the third projection function
increases while moving along a string from its initial point to its final point. The $n$-string braids
form a group, denoted by $\mathbb{B}_n$. This group is generated by $n-1$ elementary braids $\sigma_i, \ i = 1, \ldots n-1$,
each of which contains a single crossing. An important result that connects knots and braids is the Alexander's theorem \cite{kumu},
which states that all knot isotopy classes can be represented by closure of braids. Using generators of the braid group, we observe that $(\sigma_1\sigma_2\ldots \sigma_{p-1})^q$, denoted by $B(p,q)$, is a toric braid as its closure is a torus knot
$K(p,q)$. In \cite{mura}, it was proved that if $p>q$, then the closure of the braid $B(p,q)$ is a minimum crossing diagram for $K(p,q)$. In \cite{km-2, km-3}, Kronheimer and Mrowka used Gauge theory to prove that the unknotting number of an algebraic knot is equal to the genus of the Milnor fiber. A consequence of this result provides the unknotting number of a torus knot of type $(p,q)$ to be $\displaystyle ((p-1)(q-1))/2$. In general, it is not true that some minimum crossing diagram contains $u(K)$ crossings whose under/over crossing change makes it unknot. But it is interesting to observe that in $B(p,q)$, there exist $(p-1)(q-1)/2$ crossings whose change, from over to under, gives a braid whose closure is a trivial knot. A well-established technique for unknotting a knot diagram is that of converting the diagram to an ascending or descending diagram. If a knot $K$ with c crossings has a bridge with b crossings, then we can convert $K$ into an ascending or descending diagram with less than equal to $(c-b)/2$ crossing changes since switching of all the descending(ascending) crossings starting from the end(beginning) of the bridge gives a totally ascending(descending) diagram. In the case of torus knot $K(p,q)$ given by the closure of $(\sigma_{p-1}\sigma_{p-2}\ldots \sigma_1)^q$, we have $c = q(p-1)$ and $b = p-1$. Thus we can convert $K(p,q)$ into an ascending or descending diagram with less than equal to $(p-1)(q-1)/2$ crossings. Since unknotting number of $K(p,q)$ is $(p-1)(q-1)/2$, both ascending number and descending number are equal to the unknotting number. This method is very easy to understand but very time consuming in deciding the exact position of these $u(K)$ crossings as there is no symmetry/pattern in the selected crossings. Also the same procedure is not valid in unknotting torus links. In this paper, we present a method to unknot torus knots and extend the same to torus links. We also present a proof for the symmetry of $p$ and $q$ in torus knot $K(p,q)$.

In Section \ref{pre}, we introduce a few definitions and some preliminaries that we use in this paper. In Section \ref{mucdftk}, we provide a method of finding minimal unknotting crossing data for all torus knots. The main results of this section are as follows:
\begin{thm}\label{thm3}
 Let $K(p, q)$ be a torus knot with $(p, q) = 1$. Then, the following statements are equivalent:
\begin{enumerate}
 \item The unknotting number of $K(p, q)$ is equal to the number of elements in the $U-$crossing data of $B(p, q)$,
 \item $q \equiv 1 \ \textrm{or} \ p-1 \ (mod \ p)$,
 \item The $U-$crossing data of $K(p, q)$ is a minimal unknotting crossing data for $K(p, q)$.
\end{enumerate}
\end{thm}
\begin{thm}\label{thm4}  For every $p$ and $a$,\ where\ $p>a$, the\ $p-$braid
\[\eta_1\kappa_{p-1}\eta_2\kappa_{p-2}\sigma_{p-1}^{-1}\cdots \eta_a\kappa_{p-a}\sigma_{{p-a}+1}^{-1}\cdots \sigma_{p-2}^{-1}\sigma_{p-1}^{-1} \sim_{M} \eta_1\eta_2\cdots\eta_a\] where
$\eta_i=\sigma_1^{g_{i,1}}\sigma_2^{g_{i,2}}\cdots \sigma_{{p-a}-1}^{g_{i,{p-a}-1}}$,\ with $g_{i,j}=1\ or -1 $, for $i=1,2,\ldots,a$\ and $j=1,2,\ldots,{p-a}-1$, and $\kappa_{j}=\sigma_{p-a}\sigma_{{p-a}+1}\cdots \sigma_j$, for $j\geq {p-a} $.
\end{thm}
In Section \ref{pqsym}, we present a simple proof for symmetry of $p$ and $q$ in $K(p,q)$.
\begin{thm}\label{thm6} Let $B(p,a)$ be a $p$-braid, where $p>a$, then if we change crossings in $B(p,a)$, based on $U(B(p,a))$,\ we get a braid which is Markov equivalent to $B(a,p-a)$.
\end{thm}
In section \ref{mucdftl}, we extend the method of finding minimal unknotting crossing data to torus links. In Appendix A we write a matlab program(taking $p$ and $q$ as input) which gives two different sets of minimal unknotting crossing data for any torus knot $K(p,q)$.

\section{Preliminaries}\label{pre}

\begin{definition}
A crossing data for any $n$-braid $\beta_n$, denoted by $[1,2,\ldots,k]$, is a finite sequence of natural numbers enclosed in a
bracket, given to the crossings starting from the first crossing from the top to the last crossing at the bottom, based on the braid representation
of $\beta_n$ using elementary braids $\sigma_1,\sigma_2,\ldots, \sigma_{n-1}$.
\end{definition}

For example, the crossing data for the braid $\beta = \sigma_1^{-1}\sigma_2\sigma_3^{-1}(\sigma_1\sigma_2\sigma_3)^3$ is $[1,2,3,\ldots, 12]$
(Figure 1).

From the definition of crossing data for a braid, we observe the following:
\begin{enumerate}
\item Two different $n-$braids may have the same crossing data. This is mainly because the crossing data does not provide under/over crossing information at a crossing.
\item Crossing data for any two braids, having equal number of crossings, is same.
\item Crossing data for any two equivalent braids need not be the same.
\item The crossing data for $B(p, q)$ is $[1,2,\ldots, p-1,\ldots, q(p-1)].$
\end{enumerate}
\begin{figure}
\begin{center}
\includegraphics[width=3cm,height=4cm]{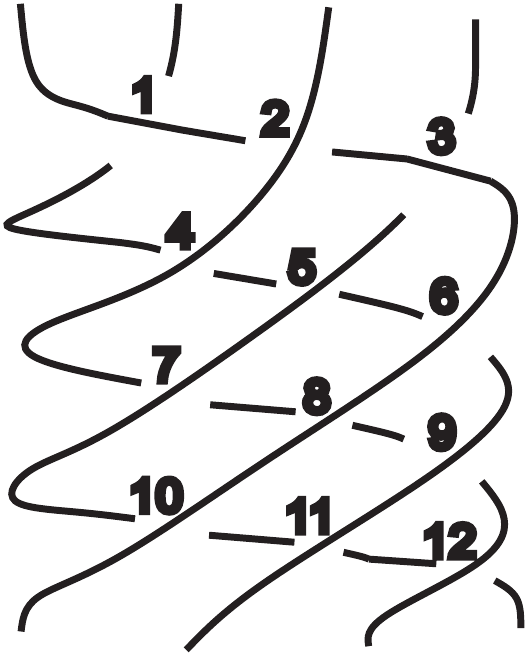}
\caption{$ \sigma_1^{-1}\sigma_2\sigma_3^{-1}(\sigma_1\sigma_2\sigma_3)^3$}
\end{center}
\end{figure}
\begin{definition}
 An unknotting crossing data for an $n-$braid, $\beta_n$, is a subsequence to the crossing data of $\beta_n$ such that if we make crossing change at these crossing positions then the closure of this braid is equivalent to unknot or unlink. If the number of elements in an unknotting crossing data
of a braid $\beta_n$ is equal to the unknotting number of $K$ (where $K$ is the closure of $\beta_n$), then this unknotting crossing data is known as minimal
unknotting crossing data for both the braid $\beta_n$ and the knot $K$.
\end{definition}

\begin{lem}\label{lem1}
For every $n$, the $(n+1)$-braid $\eta_1\eta_2\cdots\eta_n$ is equal to $\sigma_n\sigma_{n-1}\cdots\sigma_2\sigma_1$,
where $\eta_1 = \sigma_1\sigma_2 \cdots\sigma_n$, $\eta_2 = \sigma_1\sigma_2\cdots\sigma_{n-1}\sigma_n^{-1}, \ldots$,
$\eta_n = \sigma_1\sigma_2^{-1}\cdots\sigma_n^{-1}$.
\end{lem}

\begin{pf}
We prove this by means of mathematical induction. Let $P(n)$ be the statement that
$\eta_1\eta_2\cdots\eta_n$ is equal to $\sigma_n\sigma_{n-1}\cdots\sigma_2\sigma_1$.\\
Basis Step: $P(1)$ is true, because the $2$-braid $\eta_1 = \sigma_1$. Also $P(2)$ is true, because the $3$-braid
$\eta_1\eta_2 = \sigma_1\sigma_2\sigma_1\sigma_2^{-1} = \sigma_2\sigma_1\sigma_2\sigma_2^{-1} = \sigma_2\sigma_1$.
Here, we have used the fundamental relation $\sigma_1\sigma_2\sigma_1 = \sigma_2\sigma_1\sigma_2$.\\
Inductive Step: Let $K$ be some fixed positive integer and $P(k)$ is true, that is, the $(k+1)$-braid $\eta_1\eta_2\cdots\eta_k$ is equal to $\sigma_k\sigma_{k-1}\cdots\sigma_2\sigma_1$,
where $\eta_1 = \sigma_1\sigma_2\cdots\sigma_k$, $\eta_2 = \sigma_1\sigma_2\cdots\sigma_{k-1}\sigma_k^{-1}, \ldots$,
$\eta_k = \sigma_1\sigma_2^{-1}\cdots\sigma_k^{-1}$, is true.\\
To complete the inductive step, we must show that $P(k+1)$ is also true.\\
Consider $(k+2)$-braid
\[\eta_1\eta_2\cdots\eta_{k+1} = \sigma_1\sigma_2\cdots\underline{\sigma_{k+1}}\sigma_1\sigma_2\cdots\sigma_{k+1}^{-1}
\cdots\sigma_1\sigma_2^{-1}\cdots\sigma_{k+1}^{-1}.\]
Using simple fundamental relations such as $\sigma_{j}\sigma_i = \sigma_i\sigma_{j}$ for $|j-i|\geq 2$,
$\sigma_{i+1}\sigma_i\sigma_{i+1}^{-1} = \sigma_{i}^{-1}\sigma_{i+1}\sigma_{i}$, we observe that the above equation is equal to
\[\sigma_1\sigma_2\cdots\sigma_k\sigma_1\sigma_2\cdots\sigma_{k-1}\underbrace{\sigma_{k+1}\sigma_k\sigma_{k+1}^{-1}}
\cdots\sigma_1\sigma_2^{-1}\cdots\sigma_{k+1}^{-1}\]
\[=\sigma_1\sigma_2\cdots\sigma_k\sigma_1\sigma_2\cdots\sigma_{k-1}\sigma_{k}^{-1}\underline{\sigma_{k+1}\sigma_{k}}
\cdots\sigma_1\sigma_2^{-1}\cdots\sigma_{k+1}^{-1}\]
\[=\sigma_1\sigma_2\cdots\sigma_k\sigma_1\sigma_2\cdots\sigma_{k-1}\sigma_{k}^{-1}\sigma_1\cdots
\sigma_{k-2}\sigma_{k+1}\underbrace{\sigma_{k}\sigma_{k-1}\sigma_{k}^{-1}}\sigma_{k+1}^{-1}\cdots
\sigma_1\sigma_2^{-1}\cdots\sigma_{k+1}^{-1}\]
\[=\sigma_1\sigma_2\cdots\sigma_k\sigma_1\sigma_2\cdots\sigma_{k-1}\sigma_{k}^{-1}\sigma_1\cdots
\sigma_{k-2}\sigma_{k+1}\sigma_{k-1}^{-1}\sigma_{k}\sigma_{k-1}\sigma_{k+1}^{-1}\cdots
\sigma_1\sigma_2^{-1}\cdots\sigma_{k+1}^{-1}\]
\[=\sigma_1\sigma_2\cdots\sigma_k\sigma_1\sigma_2\cdots\sigma_{k-1}\sigma_{k}^{-1}\sigma_1\cdots
\sigma_{k-2}\sigma_{k-1}^{-1}\underbrace{\sigma_{k+1}\sigma_k\sigma_{k+1}^{-1}}\sigma_{k-1}\cdots
\sigma_1\sigma_2^{-1}\cdots\sigma_{k+1}^{-1} \]
\[=\sigma_1\sigma_2\cdots\sigma_k\sigma_1\sigma_2\cdots\sigma_{k-1}\sigma_{k}^{-1}\sigma_1\cdots
\sigma_{k-2}\sigma_{k-1}^{-1}\sigma_k^{-1}\underline{\sigma_{k+1}\sigma_k\sigma_{k-1}}\cdots
\sigma_1\sigma_2^{-1}\cdots\sigma_{k+1}^{-1}. \]
Finally, by symmetry, the above equation is equal to
\[\sigma_1\sigma_2\cdots\sigma_k\sigma_1\sigma_2\cdots\sigma_{k-1}\sigma_{k}^{-1}\cdots \sigma_1\sigma_2^{-1}\cdots
\sigma_k^{-1}\sigma_1^{-1}\sigma_2^{-1}\cdots\sigma_k^{-1}\underline{\sigma_{k+1}\sigma_k\cdots\sigma_1}.\]
Now, by induction hypothesis \[\sigma_1\sigma_2\cdots\sigma_k\sigma_1\sigma_2\cdots\sigma_{k-1}\sigma_{k}^{-1}\cdots \sigma_1\sigma_2^{-1}\cdots
\sigma_k^{-1} = \sigma_k\sigma_{k-1}\cdots\sigma_2\sigma_1,\]
we have
\begin{align}
\sigma_1\sigma_2\cdots\sigma_k\sigma_1\sigma_2\cdots\sigma_{k-1} &\sigma_{k}^{-1}\cdots \sigma_1\sigma_2^{-1}\cdots
\sigma_k^{-1}\sigma_1^{-1}\sigma_2^{-1}\cdots\sigma_k^{-1}\sigma_{k+1}\sigma_k\cdots\sigma_1 \notag \\
&=\sigma_k\sigma_{k-1}\cdots\sigma_2\sigma_1\sigma_1^{-1}\sigma_2^{-1}\cdots\sigma_k^{-1}\sigma_{k+1}\sigma_k\cdots\sigma_1 \notag\\
&=\sigma_{k+1}\sigma_k\cdots\sigma_1.\notag
\end{align}
Therefore, by the principle of mathematical induction, we have shown that $\eta_1\eta_2\cdots\eta_n=\sigma_n\sigma_{n-1}\cdots\sigma_2\sigma_1$
is true for all positive integers $n$.
\end{pf}

\section{A Method for Determining Minimal Unknotting Crossing Data}\label{mucdftk}

In this section, we present a method for determining the minimal unknotting crossing data for the braids of type $(\sigma_1\sigma_2\cdots\sigma_{p-1})^q$,
where $q \equiv 1 \ \textrm{or} \ (p-1) \ (mod \ p)$, so that, after switching the crossings at this selected crossing data, we obtain a braid whose closure is isotopically
equivalent to the trivial knot. Though this method does not enable us to determine the minimal unknotting crossing data for all toric braids of type $B(p, q)$, it ultimately provides us a way to find out
the minimal unknotting crossing data for $B(p, q)$.\\

The following procedure directly provides the minimal unknotting crossing data for a large class of toric braids, known as toric braids
$B(p, q)$, where $(p, q) = 1$, and $q\equiv 1 \ \textrm{or} \ p-1 \ (mod \ p)$,  and helps to find out the
minimal unknotting crossing data for every other toric braid.

\noindent \textbf{Procedure 1:} \\
Consider a toric braid of type $B(p, q)$, such that $\ (p, q) = 1$.\\
\textbf{Step 1}: First represent $B(p, q)$ as
\[ \displaystyle \underbrace{\underbrace{\sigma_1\sigma_2\cdots \sigma_{p-1}}\underbrace{\sigma_1\sigma_2\cdots \sigma_{p-1}}
\cdots \underbrace{\sigma_1\sigma_2\cdots \sigma_{p-1}}}_{q-\textrm{factors}}.\]
\textbf{Step 2}: The crossing data for the braid $B(p, q)$, based on step 1, is\\
\[ [\underbrace{1, 2, \ldots, p-1}_1,\underbrace{p, p+1, \ldots, 2(p-1)}_2, \dots \underbrace{(q-1)(p-1)+1, \ldots, q(p-1)}_q ]. \]
\textbf{Step 3}:
\begin{enumerate}
 \item When $p\geq q$: Consider the crossing data \[ [2(p-1), 3(p-1) - 1, 3(p-1), \ldots, q(p-1)-(q-2), q(p-1)-(q-3), \ldots, q(p-1)] \] if $q\neq 1$ and $[\emptyset]$ if $q=1$.
\item When $p<q$: Let $q = mp+a$, where $m\in \mathbb{N}, \ a<p$. Let
\[Y = \{2(p-1), 3(p-1) - 1, 3(p-1), \ldots, a(p-1)-(a-2), a(p-1)-(a-3), \ldots, a(p-1)\},\] if $a\neq 1$ and $\{\emptyset\}$ if $a=1$;
\[X = \{2(p-1), 3(p-1)-1, 3(p-1),\ldots, p(p-1)-(p-2), p(p-1) - (p-3), \ldots, p(p-1)\}\]
and $X+y = y +X = \{x+y\ | \ \textrm{for \ each} \ x\in X\}$. Now, consider the crossing data
\[ [X, p(p-1) + X, 2p(p-1) + X, \ldots, (m-1)p(p-1) + X, mp(p-1) + Y] .\]
\end{enumerate}
\textbf{Step 4}: Denote the crossing data obtained in Step 3 as $U-$crossing data for $B(p, q)$. Represent the $U-$crossing data for $B(p, q)$ as $U(B(p, q))$.\\

We observe that for toric braids
$B(p, q)$, where $(p, q) = 1$, and $q\equiv 1 \ \textrm{or} \ p-1 \ (mod \ p)$, the $U(B(p, q))$ is same as
minimal unknotting crossing data for $B(p, q)$ (or $K(p, q)$). Using this procedure, we can also find minimal unknotting crossing data for all torus knots of
type $K(p,q)$ with $(p, q) = 1$ and $p = 2, 3, 4,$ or $6$. In general, the $U-$crossing data for any toric braid $B(p, q)$ will help us to
determine the minimal unknotting crossing data for that braid.

\begin{thm}\label{thm1}
 For every $n$, the $(n+1)$-braid
\[\sigma_1\sigma_2\cdots\sigma_n \sigma_1\sigma_2\cdots\sigma_{n-1}\sigma_n^{-1}
\sigma_1\sigma_2\cdots\sigma_{n-1}^{-1}\sigma_n^{-1}\cdots \sigma_1^{-1}\sigma_2^{-1}\cdots\sigma_n^{-1}\]
is a trivial $(n+1)$-braid.
\end{thm}
Proof directly follows from Lemma \ref{lem1}. 

\begin{thm}\label{thm2}
 Let $K(p, q)$ be a torus knot with $(p, q) = 1$. If \[q\equiv 1 \ \textrm{or} \ p-1 \ (mod \ p),\] then the $U-$crossing data for $B(p, q)$
is a minimal unknotting crossing data for $B(p, q)$ (or $K(p, q)$).
\end{thm}

\begin{pf}
Consider $K(p, q)$.\\
\noindent\textbf{Case 1.} When $q\equiv 1 \ (mod \ p)$.\\
Then, $q = mp+1$ for some $m$.
The toric braid representation of $k(p, q)$ is $B(p, q) = \displaystyle \underbrace{\underbrace{\sigma_1\sigma_2\cdots \sigma_{p-1}}\underbrace{\sigma_1\sigma_2\cdots \sigma_{p-1}}
\cdots \underbrace{\sigma_1\sigma_2\cdots \sigma_{p-1}}}_{q-\textrm{factors}}.$
By making crossing changes at the $U-$crossing data for $B(p, q)$, we obtain:
\begin{align}
&\displaystyle \underbrace{\underbrace{\sigma_1\sigma_2\cdots \sigma_{p-1}}\underbrace{\sigma_1\sigma_2\cdots \sigma_{p-1}^{-1}}
\cdots \underbrace{\sigma_1^{-1}\sigma_2^{-1}\cdots \sigma_{p-1}^{-1}}\underbrace{\sigma_1\sigma_2\cdots
\sigma_{p-1}}}_{q = (mp+1)-\textrm{factors}} \notag \\
&=\displaystyle{(\underbrace{\sigma_1\sigma_2\cdots \sigma_{p-1}}\underbrace{\sigma_1\sigma_2\cdots \sigma_{p-1}^{-1}}\cdots \underbrace{\sigma_1^{-1}\sigma_2^{-1}\cdots \sigma_{p-1}^{-1}})}^m \underbrace{\sigma_1\sigma_2\cdots
\sigma_{p-1}} \notag \\
&=\underbrace{\sigma_1\sigma_2\cdots\sigma_{p-1}}~~(by\ Theorem\ \ref{thm1})\notag
\end{align}
Now, we show that the number of elements in the
$U-$crossing data is equal to the unknotting number for these torus knots.
From \cite{km-2, km-3}, we know that the unknotting number for these torus knots is equal to
\[\displaystyle \frac{(p-1)(q-1)}{2} = \frac{(p-1)mp}{2}.\]
The number of elements in the $U-$crossing data is equal to \[\displaystyle m \sum_{i=0}^{p-1} i = \frac{m\cdot (p-1)p}{2} = \frac{(p-1)mp}{2},\]
which is same as the unknotting number. \\

\noindent \textbf{Case 2.} When $q \equiv p-1 \ (mod \ p)$.\\
Then, $q = mp-1$ for some $m$.
The toric braid representation of $k(p, q)$ is $B(p, q) = (\sigma_1\sigma_2\cdots \sigma_{p-1})^{mp-1}$.
By using procedure 1, making the crossing changes at the $U-$crossing data, we obtain:
\[\displaystyle \underbrace{\underbrace{\sigma_1\sigma_2\cdots \sigma_{p-1}}\underbrace{\sigma_1\sigma_2\cdots \sigma_{p-1}^{-1}}
\cdots \underbrace{\sigma_1^{-1}\sigma_2^{-1}\cdots \sigma_{p-1}^{-1}}\underbrace{\sigma_1\sigma_2\cdots \sigma_{p-1}}\cdots
\underbrace{\sigma_1\sigma_2^{-1}\cdots \sigma_{p-1}^{-1}}}_{q = (mp-1)-\textrm{factors}}\]
$= \alpha^{m-1} \beta$, where:
\[\alpha = \sigma_1\sigma_2\cdots \sigma_{p-1}\sigma_1\sigma_2\cdots \sigma_{p-1}^{-1}
\sigma_1\sigma_2\cdots\sigma_{p-2}^{-1}\sigma_{p-1}^{-1}\cdots \sigma_1^{-1}\sigma_2^{-1}\cdots
\sigma_{p-1}^{-1}\]
\[\beta = \sigma_1\sigma_2\cdots \sigma_{p-1}\sigma_1\sigma_2\cdots \sigma_{p-1}^{-1}\sigma_1\sigma_2\cdots \sigma_{p-2}^{-1}\sigma_{p-1}^{-1}
\cdots \sigma_1\sigma_2^{-1}\cdots \sigma_{p-1}^{-1} .\]
By Theorem \ref{thm1}, we get
\[\alpha^{m-1} \beta =\beta= \sigma_1\sigma_2\cdots \sigma_{p-1}\sigma_1\sigma_2\cdots \sigma_{p-1}^{-1}\sigma_1\sigma_2\cdots \sigma_{p-2}^{-1}\sigma_{p-1}^{-1}
\cdots \sigma_1\sigma_2^{-1}\cdots \sigma_{p-1}^{-1}.\]
By Lemma \ref{lem1}, the above equation is equal to $\sigma_{p-1}\sigma_{p-2}\cdots \sigma_2 \sigma_1$, whose closure is a trivial knot.\\
Now, we show that the number of elements in the $U-$crossing data is equal to the unknotting number for these torus knots. Observe that the unknotting number is equal to $\displaystyle \frac{(p-1)(mp-2)}{2}$ and the number of elements in the $U-$crossing data is equal to
\[\displaystyle (m-1) \sum_{i=0}^{p-1} i + \sum_{i=0}^{p-2} i = \frac{(m-1)\cdot (p-1)p}{2} + \frac{(p-2)(p-1)}{2} = \frac{(p-1)(mp-2)}{2}.\]
\end{pf} \vspace{0.3cm}
Now, we have sufficient data to prove the main theorem of this section.
\newproof{pot1}{Proof of Theorem \ref{thm3}}
\begin{pot1} $(1) \rightarrow (2)$: For given torus knot $K(p, q)$, we can write $q = mp + a$ for some integer $m$ and for some $0<a<p$. Observe that the unknotting number of $K(p, q)$ is equal to $\displaystyle \frac{mp(p-1)}{2} + \frac{(p-1)(a-1)}{2}$ and the number of elements in the $U-$crossing data is
$\displaystyle m\cdot \sum_{i=0}^{p-1} i + \sum_{i=0}^{a-1} i = \frac{mp(p-1)}{2} + \frac{a(a-1)}{2}$.
By hypothesis, these two are equal. Therefore, we have two solutions for $a (< p)$, i.e., either
$a=1$ or $a = p-1$. Hence, $q \equiv 1 \ \textrm{or} \ p-1 \ (mod \ p)$.
\noindent $(2) \rightarrow (3)$ is nothing but Theorem \ref{thm2} and $(3) \rightarrow (1)$ directly follows from the definition of
 minimal unknotting crossing data.
\end{pot1}

\begin{definition}
The braid $\overline{\beta}$ is the braid obtained from $\beta$ by rotating $\beta$ in $\mathbb{R}^3$
about the $y$-axis through an angle of $\pi$.
\end{definition}
 \begin{rmk}\label{rmk2}
(a). It is easy to observe that the closure of the braids $\beta$ and $\overline{\beta}$ are equivalent. Actually  $\beta$ and  $\overline{\beta}$ are same because we obtain  $\overline{\beta}$ from  $\beta$ without any braid operation. However we can not say that for any braid $\eta $, $\eta\beta\sim_{M}\eta\overline{\beta}$.  Also, observe that if
$\beta = B(p, q)= (\sigma_1\sigma_2\cdots \sigma_{p-1})^q$, then
$\overline{\beta} = \overline{B}(p,q) = (\sigma_{p-1}\sigma_{p-2}\cdots \sigma_1)^q$.

 (b). If $X$ is the $U-$crossing data of $B(p,q)$, then the elements in $U-$crossing data of $\overline{B}(p,q)$ are $(p-1)q + 1-X$.
\end{rmk}

\begin{lem}\label{lem2} If $1\leq j<i\leq n$, then
\[\sigma_i^g\sigma_j^{g_j}\sigma_{j+1}^{g_{j+1}}\cdots\sigma_{i-1}^{g_{i-1}}\sigma_i^{g_i}\sigma_{i+1}^{g_{i+1}}
\cdots\sigma_n^{g_n}=\sigma_j^{g_j}\sigma_{j+1}^{g_{j+1}}\cdots\sigma_{i-1}^{g_i}\sigma_i^{g_{i-1}}
\sigma_{i+1}^{g_{i+1}}\cdots\sigma_n^{g_n}\sigma_{i-1}^{g}\]
holds if either $g=g_{i-1}$ or $g_{i-1}=g_{i}$, where $g, g_k\ (k=j,j+1,\ldots,n)$ have values $1\ or\ -1$.
\end{lem}

\begin{pf} Proof directly follows by using the braid relations
\[\sigma_j\sigma_i=\sigma_i\sigma_j\ if\ |i-j|>1\ and\ \sigma_{i+1}^{\alpha_{1}}\sigma_i^{\alpha_2}\sigma_{i+1}^{\alpha_{3}}=\sigma_i^{\alpha_{3}}\sigma_{i+1}^{\alpha
_2}\sigma_i^{\alpha_{1}}\]
if either $\alpha_1=\alpha_2$ or $\alpha_2=\alpha_3$, where $\alpha_1,\alpha_2\ and\ \alpha_3\ are\ 1\ or\ -1.$\\
We can verify the last relation for different possible values of $\alpha_1,\alpha_2\ and\ \alpha_3$.
\end{pf} \vspace{0.3cm}
Now, we prove second main theorem of this section.
\newproof{pot2}{Proof of Theorem \ref{thm4}}
\begin{pot2} Consider, the $p-$braid \[\eta_1\kappa_{p-1}\eta_2\kappa_{p-2}\sigma_{p-1}^{-1}\cdots \eta_a\kappa_{p-a}\sigma_{{p-a}+1}^{-1}\cdots \sigma_{p-2}^{-1}\sigma_{p-1}^{-1}\]

$\sim_{M}\eta_1\kappa_{p-2}\eta_2\kappa_{p-3}\sigma_{p-2}^{-1}\underline{\sigma_{p-1}\sigma_{p-2}}\cdots \eta_a\kappa_{p-a}\sigma_{{p-a}+1}^{-1}\cdots \sigma_{p-2}^{-1}\sigma_{p-1}^{-1}$

$\sim_{M}\eta_1\kappa_{p-2}\eta_2\kappa_{p-3}\sigma_{p-2}^{-1}\eta_3\kappa_{p-4}\sigma_{p-3}^{-1}\sigma_{p-2}^{-1}
\underline{\sigma_{p-1}\sigma_{p-2}\sigma_{p-3}}\cdots \eta_a\kappa_{p-a}\sigma_{{p-a}+1}^{-1}\cdots \sigma_{p-2}^{-1}\sigma_{p-1}^{-1}$
\begin{center}
$\vdots$
\end{center}

$\sim_{M}\eta_1\kappa_{p-2}\eta_2\kappa_{p-3}\sigma_{p-2}^{-1}\eta_3\kappa_{p-4}\sigma_{p-3}^{-1}\sigma_{p-2}^{-1}
\cdots \eta_a\sigma_{p-a}^{-1}\sigma_{{p-a}+1}^{-1}\cdots \sigma_{p-2}^{-1}\underline{\sigma_{p-1}\sigma_{p-2}\cdots\sigma_{p-a}}$

$\sim_{M}\eta_1\kappa_{p-2}\eta_2\kappa_{p-3}\sigma_{p-2}^{-1}\eta_3\kappa_{p-4}\sigma_{p-3}^{-1}\sigma_{p-2}^{-1}
\cdots \eta_{a-1}\sigma_{p-a}\sigma_{{p-a}+1}^{-1}\cdots \sigma_{p-2}^{-1}\eta_{a}$.\\
Similarly, we obtain

$\sim_{M}\eta_1\kappa_{p-3}\eta_2\kappa_{p-4}\sigma_{p-3}^{-1}\eta_3\kappa_{p-5}\sigma_{p-4}^{-1}\sigma_{p-3}^{-1}
\cdots \eta_{a-2}\sigma_{p-a}\sigma_{{p-a}+1}^{-1}\cdots \sigma_{p-3}^{-1}\eta_{a-1}\eta_{a}$.\\
Finally, we have

$\sim_{M}\eta_1\eta_2\cdots\eta_a$.\\

\end{pot2}

\begin{rmk}\label{rmk3}   If we take $g_{i,j}=1 $ for $i=1,2,\ldots,a$\ and $j=1,2,\ldots,n-1$\, we obtain
$\sigma_1\sigma_2\cdots \sigma_{p-2}\sigma_{p-1}\sigma_1\sigma_2\cdots \sigma_{p-2}\sigma_{p-1}^{-1}\cdots \sigma_1\sigma_2\cdots\sigma_{p-a}\sigma_{p-a+1}^{-1}\cdots \sigma_{p-2}^{-1}\sigma_{p-1}^{-1}$\\
$\sim_{M}\sigma_1\sigma_2\cdots\sigma_{p-a-1}\sigma_1\sigma_2\cdots \sigma_{p-a-1}\cdots \sigma_1\sigma_2\cdots\sigma_{p-a-1}\ =B(p-a,a)$

That is, if we change crossings in $B(p,a)$, based on $U(B(p,a))$, we obtain a diagram which is Markov equivalent to $B(p-a,a).$
\end{rmk}
\begin{rmk}\label{rmk4} For every $p,q$, with $(p,q)=1 \ and\ q \equiv a (mod\ p)$,
\[B(p,q)= {(\sigma_1\sigma_2 \cdots \sigma_{p-2}\sigma_{p-1})}^{mp+a};for\ some\ m\geq 0, a<p.\]
Then, by changing the crossings in $B(p,q)$, based on $U(B(p,q))$, we get\\
$\alpha ^m \beta $, where
$\alpha=\sigma_1\sigma_2 \cdots \sigma_{p-2}\sigma_{p-1}\sigma_1\sigma_2 \cdots \sigma_{p-2}\sigma_{p-1}^{-1}\cdots \sigma_1^{-1}\sigma_2^{-1} \cdots \sigma_{p-2}^{-1}\sigma_{p-1}^{-1}$
and  $\beta = {\sigma_1\sigma_2 \cdots \sigma_{p-2}\sigma_{p-1}\sigma_1\sigma_2 \cdots \sigma_{p-2}\sigma_{p-1}^{-1}\cdots \sigma_1\sigma_2 \cdots\sigma_{p-(a-1)}^{-1}\sigma_{p-(a-2)}^{-1}\cdots \sigma_{p-2}^{-1}\sigma_{p-1}^{-1}}$.\\
By Theorem \ref{thm1} and Remark \ref{rmk3}, $\alpha $ is a trivial $p-$braid and $\beta$ is Markov equivalent to $B(p-a,a)$.
Thus, we can find a minimal unknotting crossing data for $K(p,q)$, if we can find the unknotting crossing data for
$B(p-a,a)$ such that the number of elements in this unknotting crossing data is equal to the unknotting number of $K(p, q) \  \smallsetminus  |U(B(p,q))|$, which is equal to
\[ \frac{(p-a-1)(a-1)}{2}.\]
Observe that, this is the unknotting number of $K(p-a,a)$.
\end{rmk}

Now we will provide the unknotting process. Consider any torus knot of type $K(p, q)$, where $(p, q) = 1\ and\ q=mp+a$. Let $p_1=p$ and $q_1=q$. Since $(p_1,q_1)=1$, Euclid's algorithm ensures that $\exists\ n\in N$ such that if
we define $q_i,p_i$ as
\[q_{i+1}=q_i-s_i p_i;\ q_{i+1}<p_i\]
\[p_{i+1}=p_i-m_i q_{i+1};\ p_{i+1}<q_{i+1}\]
we get $ q_n\equiv 1\ or\ (p_n-1)\ (mod\ p_n )$.\\
Now, if we select $U$-crossing data for $K(p_i,q_i), K(p_i-q_{i+1},q_{i+1}), K(p_i-2q_{i+1},q_{i+1}),\cdots, K(p_i-(m_i-1)q_{i+1},q_{i+1})$ for $i=1,2,\cdots,n$, then by Theorem \ref{thm4} and Remark \ref{rmk4}, there are corresponding crossings in $K(p,q)$ which provides minimal unknotting crossing data for $K(p,q)$. These corresponding crossings are determined as below:

For finding the minimal unknotting crossing data for a torus knot we consider the following steps.

Let $K(p,q)$\ be a torus knot. Then, (say)

 \textbf{Step 1}: $(p_1,q_1)=(p,q)\ ;  \ where\ q_1=m_1p_1+a_1$

 \textbf{Step 2}: $(p_2,q_2)=(p,a_1)\ ;  \ where\  p_2=m_2q_2+a_2$

 \textbf{Step 3}: $(p_3,q_3)=(a_2,a_1)\ ;  \ where\ q_3=m_3p_3+a_3$

 \textbf{Step 4}: $(p_4,q_4)=(a_2,a_3)\ ; \ where\  p_4=m_4q_4+a_4$
\begin{center}
$ \vdots $
\end{center}

\textbf{Step i}: $
 (p_i,q_i)=
 \begin{cases}
 (a_{i-1},a_{i-2}); \  where\ i \ is\ odd\ and\ q_i=m_ip_i+a_i \\
 (a_{i-2},a_{i-1}); \ where\ if\ i\ is\ even\ and\ p_i=m_iq_i+a_i 

 \end{cases}
 $
\begin{center}
$ \vdots $
\end{center}

\textbf{Step n}: $
 (p_n,q_n)=
 \begin{cases}
 (a_{n-1},a_{n-2}); \  where\ if\ n\ is\ odd\ and\ q_n\equiv 1\ or\ (p_n-1)\ (mod\ p_n )\\
 (a_{n-2},a_{n-1});\ where\ if\ n\ is\ even\ and\ p_n\equiv 1\ (mod\ q_n ) 

 \end{cases}
 $\\
is our last step.\\
Now, consider $B_{i=odd}=\{j p_i(p-1)+X_i : j=0,1,\ldots,m_i\}$; and $B_{i=even}=\{X_{i,j}: j=1,2,\ldots ,m_i-1\}$. Then the minimal unknotting crossing data for $K(p,q)$ is  \[[B_1,m_1p_1(p-1)+\{B_2,B_3,m_3p_3(p-1)+\{B_4,B_5,\ldots  , m_{(n-2)}p_{(n-2)}(p-1)+\{B_{(n-1)},B_n\}\}\ldots \}]\] if n is odd; and \[[B_1,m_1p_1(p-1)+\{B_2,B_3,m_3p_3(p-1)+\{B_4,B_5,\ldots  , m_{(n-1)}p_{(n-1)}(p-1)+\{B_n\}\}\ldots \}]\] if n is even.

Here $X_i=\{(p-1)+p_i-1,2(p-1)+p_i-1,2(p-1)+p_i-2,\ldots, (l-1)(p-1)+p_i-1,(l-1)(p-1)+p_i-2,\ldots ,(l-1)(p-1)+p_i-(l-1)\}$,
where $l=p_i\ for\ j=0,1,\ldots ,m_i-1;\ and\ l=a_i\ for\ j=m_i$;\\ and
$X_{i,j}=\{(p-1)+(p_i-jq_i)-1,2(p-1)+(p_i-jq_i)-1,2(p-1)+(p_i-jq_i)-2,\ldots, (q_i-1)(p-1)+(p_i-jq_i)-1,(q_i-1)(p-1)+(p_i-jq_i)-2,\ldots ,(q_i-1)(p-1)+(p_i-jq_i)-(q_i-1)\}.$

\begin{figure}
    \centering
   \includegraphics[height=12.6cm,width=5.8in]{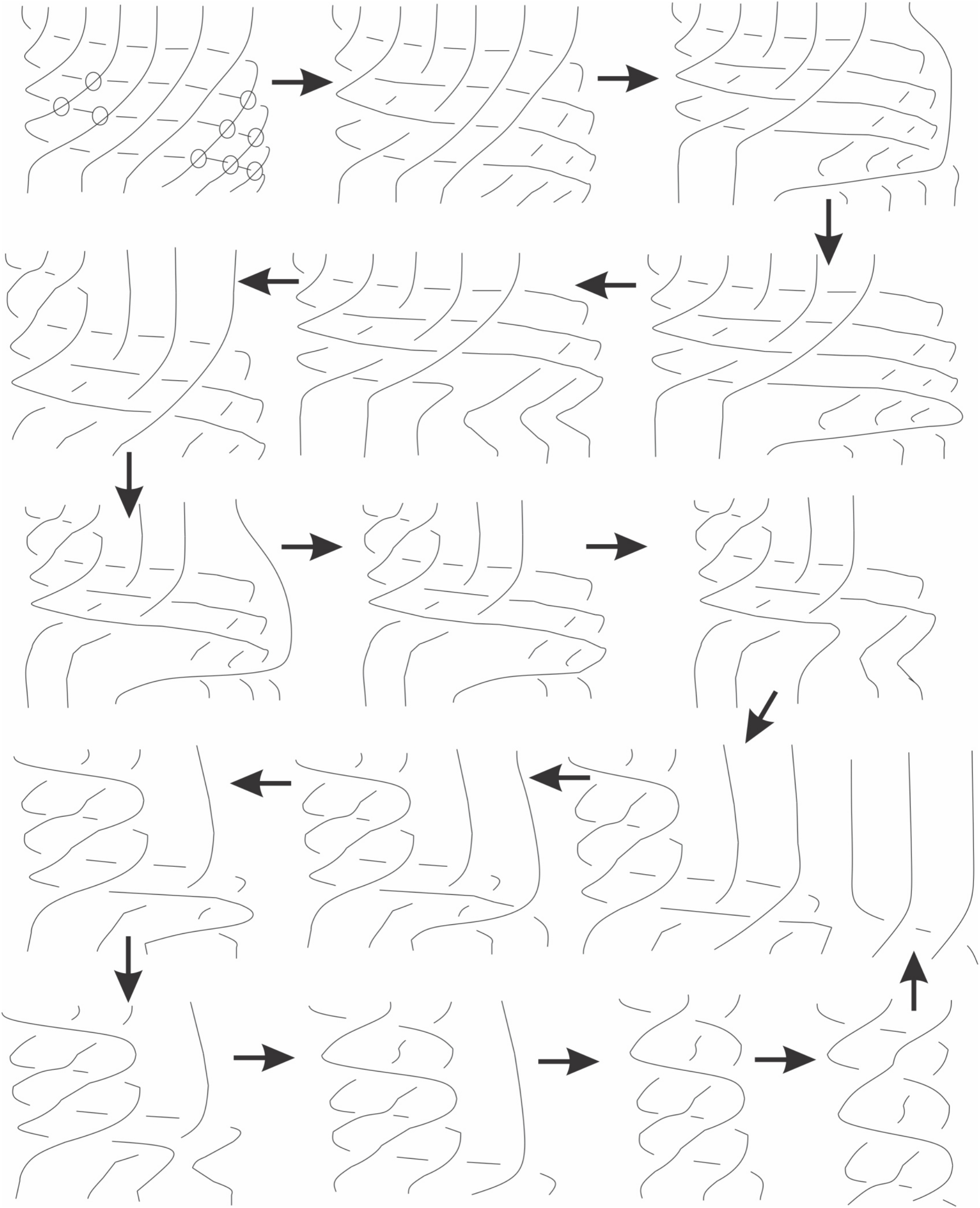}
    \caption{Unknotting procedure for torus knot $K(7,4)$  }
    \label{FF}
\end{figure}

\begin{ex}\label{ex2} To find the unknotting procedure and minimal unknotting crossing data for torus knot $K(7,4)$.\\

By the above described procedure, we have

1. $(p_1,q_1)=(p,q)=(7,4)\ ; $ \ where $4=0.7+4$

 2. $(p_2,q_2)=(p,a_1)=(7,4)\ ; $ \ where $ 7=1.4+3$

 3. $(p_n,q_n)=(p_3,q_3)=(3,4)\ ; $ \ where $ 4=1.3+1$\\
Here $p=p_1=p_2=7; q=q_1=q_2=q_3=a_1=4;p_3=a_2=3;m_1=0;m_2=m_3=1\ and\ a_3=1.$\\
For $n=3$, the minimal unknotting crossing data is $[B_1,m_1p_1(p-1)+\{B_2,B_3\}]=[8,12,13,14,17,18,22,23,24]$,
where $B_1=\{12,17,18,22,23,24 \}, B_2=\{\emptyset\}$\ and $B_3=\{8,13,14\}$.
\end{ex}
\begin{ex}\label{ex3} Minimal unknotting crossing data for torus knot $K(13,3)$ is $[15,18,21,24,26,27,29,30,32,33,35,36]$. In this case $B_1=\{24,35,36 \}\ and\ B_2=\{21,32,33,18,29,30,15,26,27\}$.
\end{ex}
\section{{Symmetry of p\ and\ q in a torus knot}}\label{pqsym}
In this section we show that $K(p,q)=K(q,p)$, by showing that both the braids $B(p,q)$ and $B(q,p)$ are Markov equivalent to a braid $\beta$, which is obtained by changing the crossings in $B(p+q,q)$, based on $U(B(p+q,q))$.

\begin{thm}\label{thm5}
For any two natural numbers $a$ and $n$,
\[\displaystyle B(a,a) * \overline{B}(a,n)  \sim_M B(a+n,a) \]
i.e., The braid product of two toric braids $B(a,a)$ and $\overline{B}(a,n)$ is Markov equivalent to the toric braid $B(a+n,a)$.
\end{thm}
\begin{pf}
To prove this theorem, we use Markov moves and the following simple straight forward results:

\begin{enumerate}
 \item For any $i \geq 2$, $\sigma_i\sigma_1\sigma_2\cdots \sigma_n = \sigma_1\sigma_2\cdots\sigma_n\sigma_{i-1}$.
 \item If $1 \leq j < i \leq n$, then $\sigma_i\sigma_j\sigma_{j+1} \cdots \sigma_n = \sigma_j\sigma_{j+1}\cdots \sigma_n\sigma_{i-1}$.
 \end{enumerate}
It is easy to prove the above two results by using the fundamental braid relations
\[\sigma_j\sigma_i = \sigma_i\sigma_j \ \textrm{if} \ |i-j|>1, \ \textrm{and}\  \sigma_i\sigma_{i+1}\sigma_i = \sigma_{i+1}
\sigma_i\sigma_{i+1}.\]
Observe that for $n = 1$, we have
\[B(a,a) * \overline{B}(a,1) = (\sigma_1\sigma_2\cdots\sigma_{a-1})^a\sigma_{a-1}\sigma_{a-2}\cdots\sigma_1.\]
By $M_2$, we have \[B(a,a) * \overline{B}(a,1) \sim_M  (\sigma_1\sigma_2\cdots\sigma_{a-1})^a \underline{\sigma_a} \sigma_{a-1}\sigma_{a-2}\cdots\sigma_1
= (\sigma_1\sigma_2\cdots\sigma_a)^a = B(a+1, a.)\]
Similarly, for $n=2$, we have
\[B(a,a) * \overline{B}(a,2) = (\sigma_1\sigma_2\cdots\sigma_{a-1})^a\sigma_{a-1}\sigma_{a-2}\cdots\sigma_1\sigma_{a-1}\sigma_{a-2}\cdots\sigma_1.\]
By $M_2$, we have \[B(a,a) * \overline{B}(a,2) \sim_M  (\sigma_1\sigma_2\cdots\sigma_{a-1})^a \sigma_{a-1}\sigma_{a-2}
\cdots\sigma_1\underline{\sigma_a}\sigma_{a-1}\sigma_{a-2} \cdots\sigma_1\]
\[\sim_M (\sigma_1\sigma_2\cdots\sigma_{a-1})^a \sigma_a\sigma_{a-1}\sigma_{a-2}
\cdots\sigma_1\underline{\sigma_{a+1}}\sigma_a\sigma_{a-1}\sigma_{a-2} \cdots\sigma_2\]
\[\sim (\sigma_1\sigma_2\cdots\sigma_a)^a \sigma_{a+1}\sigma_a\sigma_{a-1}\sigma_{a-2} \cdots\sigma_2
= (\sigma_1\sigma_2\cdots\sigma_{a+1})^a = B(a+2, a).\]
Now, in general:
\[\displaystyle B(a,a) * \overline{B}(a,n) = (\sigma_1\sigma_2\cdots\sigma_{a-1})^a\underbrace{\underbrace{\sigma_{a-1}\sigma_{a-2}\cdots\sigma_1}\cdots
\underbrace{\sigma_{a-1}\sigma_{a-2}\cdots\sigma_1}}_{n-factors}\]
\[\displaystyle \sim_M(\sigma_1\sigma_2\cdots\sigma_{a-1})^a\underbrace{\sigma_{a-1}\sigma_{a-2}\cdots\sigma_1}\cdots
\underbrace{\sigma_{a-1}\sigma_{a-2}\cdots\sigma_1}\underline{\sigma_a}\underbrace{\sigma_{a-1}\sigma_{a-2}\cdots\sigma_1}\]
\[\displaystyle \sim_M(\sigma_1\sigma_2\cdots\sigma_{a-1})^a\underbrace{\sigma_a\sigma_{a-1}\sigma_{a-2}\cdots\sigma_1}
\underbrace{\sigma_{a}\sigma_{a-1}\cdots\sigma_2}\cdots\underbrace{\sigma_{a}\sigma_{a-1}\cdots\sigma_2}\underline{\sigma_{a+1}}\underbrace{\sigma_a\sigma_{a-1}\cdots\sigma_2}\]
\[\displaystyle \sim_M(\sigma_1\sigma_2\cdots\sigma_{a})^a\underbrace{\sigma_{a+1}\sigma_a\cdots\sigma_2}\underbrace{\sigma_{a+1}\sigma_a\cdots\sigma_3}\cdots
\underbrace{\sigma_{a+1}\sigma_a\cdots\sigma_3}\underline{\sigma_{a+2}}\underbrace{\sigma_{a+1}\sigma_a\cdots\sigma_3}\]
\hspace{2cm}$\vdots$ \hspace{1cm} $\vdots$ \hspace{1cm} $\vdots$ \hspace{1cm} $\vdots$ \hspace{1cm} $\vdots$ \hspace{1cm} $\vdots$
\hspace{1cm} $\vdots$ \hspace{1cm} $\vdots$ \hspace{1cm} $\vdots$\\
$\displaystyle \sim_M (\sigma_1\sigma_2\cdots\sigma_{a+n-2})^a\sigma_{a+n-1}\underbrace{\sigma_{a+n-2}\cdots\sigma_n}$\\
$\displaystyle \sim_M (\sigma_1\sigma_2\cdots\sigma_{a+n-1})^a$
$ = B(a+n, a)$.
\end{pf}

\begin{rmk}\label{rmk5} For $a+n=p$, by Theorem \ref{thm5}, the $p$-braid $B(p,a)$, where $p>a$, is Markov equivalent to a $a$-braid $B(a,a) * \overline{B}(a,p-a)$.
\end{rmk}

\begin{ex}\label{ex1}
Consider a toric braid of type $B(5,4)$. Then, by Theorem \ref{thm5},  $B(4,4) * \overline{B}(4,1) \sim_M B(5,4)$.
Here, in this example, we show $B(5,4) \sim_M B(4,4) * \overline{B}(4,1)$. Observe that
\[B(5,4) \sim_M \sigma_1\sigma_2\sigma_3\sigma_4\sigma_1\sigma_2\sigma_3\sigma_4\sigma_1\sigma_2\sigma_3\sigma_4\sigma_1\sigma_2\sigma_3\sigma_4\]
The same is shown in Figure 4 $(a)$. Using the equality $\sigma_4\sigma_1\sigma_2\sigma_3\sigma_4 = \sigma_1\sigma_2\sigma_3\sigma_4\sigma_3$ in the above
equation, we get \\
$B(5, 4) \sim \sigma_1\sigma_2\sigma_3\sigma_1\sigma_2\sigma_3\sigma_4\sigma_3\sigma_1\sigma_2\sigma_3\sigma_4\sigma_1\sigma_2\sigma_3\sigma_4$. This is
shown in Figure 4 $(b)$. Again, using $\sigma_4\sigma_3\sigma_1\sigma_2\sigma_3\sigma_4 = \sigma_1\sigma_2\sigma_3\sigma_4\sigma_3\sigma_2$,
we get\\
$B(5,4)\sim \sigma_1\sigma_2\sigma_3\sigma_1\sigma_2\sigma_3\sigma_1\sigma_2\sigma_3\sigma_4\sigma_3\sigma_2\sigma_1\sigma_2\sigma_3\sigma_4$.
This expression is shown in Figure 4 $(c)$. Figure 4 $(d)$ represents the following expression, which is also obtained by a similar argument,\\
\[B(5,4) \sim \sigma_1\sigma_2\sigma_3\sigma_1\sigma_2\sigma_3\sigma_1\sigma_2\sigma_3\sigma_1\sigma_2\sigma_3\sigma_4\sigma_3\sigma_2\sigma_1.\]

Now using the Markov move $M_1$, $\beta \sim_M \gamma \beta \gamma^{-1}$, where $\beta = (\sigma_1\sigma_2\sigma_3)^4\sigma_4\sigma_3\sigma_2\sigma_1$
and $\gamma = \sigma_3\sigma_2\sigma_1$, we observe that $B(5,4) \sim_M \sigma_3\sigma_2\sigma_1\sigma_1\sigma_2\sigma_3\sigma_1\sigma_2\sigma_3\sigma_1
\sigma_2\sigma_3\sigma_4$. This is shown in Figure 4 $(e)$.\\
Using Markov move $M_2$, we remove $\sigma_4$ from the above expression, so that we obtain
$B(5,4) \sim_M \sigma_3\sigma_2\sigma_1\sigma_1\sigma_2\sigma_3\sigma_1\sigma_2\sigma_3\sigma_1\sigma_2\sigma_3$, which is shown in Figure 4 $(f)$.\\
Finally, again by using Markov move $M_1$, we obtain that\\
$B(5,4)\sim_M \sigma_1\sigma_2\sigma_3\sigma_1\sigma_2\sigma_3\sigma_1\sigma_2\sigma_3\sigma_3
\sigma_2\sigma_1 = B(4,4) * \overline{B}(4,1)$. This is shown in Figure 4 $(g)$.
\end{ex}
\begin{figure}
  \begin{center}
   \includegraphics[height=8cm,width=14.5cm]{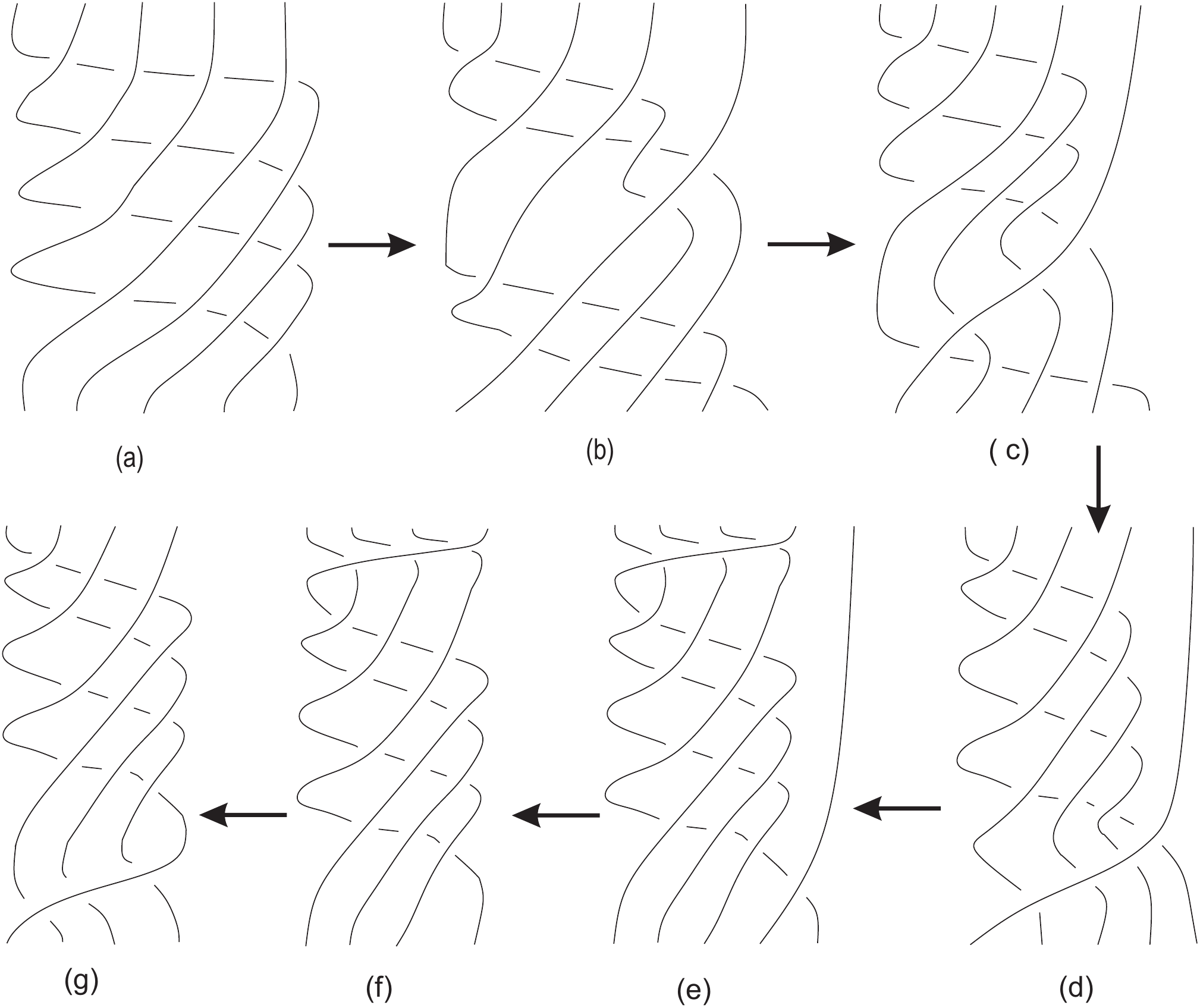}
   \caption{$B(5,4)=B(4,4)*\overline{B}(4,1)$}
 \end{center}
\end{figure} \vspace{0.3cm}
Now we prove the main theorem of this section.
\newproof{pot3}{Proof of Theorem \ref{thm6}}
\begin{pot3} Let $\beta$ is the braid obtained from $B(p,a)$ after changing the crossings in $B(p,a)$, based on $U(B(p,a))$.
Then\\
$\beta = {\sigma_1\sigma_2 \cdots \sigma_{p-2}\sigma_{p-1}\sigma_1\sigma_2 \cdots \sigma_{p-2}\sigma_{p-1}^{-1}\cdots \sigma_1\sigma_2 \cdots\sigma_{p-(a-1)}^{-1}\sigma_{p-(a-2)}^{-1}\cdots \sigma_{p-2}^{-1}\sigma_{p-1}^{-1}}$
$= \underbrace{\sigma_1\sigma_2 \cdots \sigma_{a-1}}\underline{\sigma_{a}\sigma_{a+1}\cdots\sigma_{p-1}}\underbrace{\sigma_1\sigma_2 \cdots \sigma_{p-2}\sigma_{p-1}^{-1}}\cdots \underbrace{\sigma_1\sigma_2 \cdots\sigma_{p-(a-1)}^{-1}\sigma_{p-(a-2)}^{-1}\cdots \sigma_{p-2}^{-1}\sigma_{p-1}^{-1}}$
$\sim_{M}\underbrace{\sigma_1\sigma_2 \cdots \sigma_{a-1}}\underbrace{\sigma_1\sigma_2 \cdots \sigma_{a-1}^{-1}}\underline{\sigma_{a}\sigma_{a+1}\cdots\sigma_{p-1}\sigma_{a-1}\sigma_{a}\cdots \sigma_{p-2}}\cdots
\underbrace{\sigma_1\sigma_2 \cdots\sigma_{p-(a-1)}^{-1}\cdots \sigma_{p-1}^{-1}}$
\begin{center}
$\vdots$
\end{center}
$\sim_{M}\underbrace{\sigma_1\sigma_2 \cdots \sigma_{a-1}}\underbrace{\sigma_1\sigma_2 \cdots \sigma_{a-1}^{-1}}\cdots \underbrace{\sigma_1^{-1}\sigma_2^{-1} \cdots \sigma_{a-1}^{-1}}\underline{\sigma_{a}\sigma_{a+1}\cdots\sigma_{p-1}\sigma_{a-1}\sigma_{a}\cdots \sigma_{p-2}\cdots}$\\$\underline{\cdots \sigma_{1}\sigma_{2}\cdots \sigma_{p-a}}.$\\
By Theorem \ref{thm1}\\
$\sim_{M}\underbrace{\underbrace{\sigma_{a}\sigma_{a+1}\cdots\sigma_{p-1}}\underbrace{\sigma_{a-1}\sigma_{a}\cdots \sigma_{p-2}}\cdots \underbrace{\sigma_{1}\sigma_{2}\cdots \sigma_{p-a}}}_{a-factors}$\\
$\sim_{M}\underbrace{\underbrace{\sigma_{a}\sigma_{a-1}\cdots\sigma_{1}}\underbrace{\sigma_{a+1}\sigma_{a}\cdots \sigma_{2}} \cdots\underbrace{\sigma_{p-2}\sigma_{p-3}\cdots \sigma_{p-a-1}} \underbrace{\sigma_{p-1}\sigma_{p-2}\cdots \sigma_{p-a}}}_{(p-a)-factors}$\\
$\sim_{M}\underbrace{\sigma_{a}\sigma_{a-1}\cdots\sigma_{1}}\underbrace{\sigma_{a+1}\sigma_{a}\cdots \sigma_{2}} \cdots\underbrace{\sigma_{p-2}\sigma_{p-3}\cdots \sigma_{p-a-1}} \underline{\sigma_{p-1}}\underbrace{\sigma_{p-2}\cdots \sigma_{p-a}}$\\
$\sim_{M}\underbrace{\sigma_{a}\sigma_{a-1}\cdots\sigma_{1}}\underbrace{\sigma_{a+1}\sigma_{a}\cdots \sigma_{2}} \cdots\underbrace{\sigma_{p-2}\sigma_{p-3}\cdots \sigma_{p-a-1}} \underbrace{\sigma_{p-2}\cdots \sigma_{p-a}}$\\
$\sim_{M}\underbrace{\sigma_{a}\sigma_{a-1}\cdots\sigma_{1}}\underbrace{\sigma_{a+1}\sigma_{a}\cdots \sigma_{2}} \cdots \underbrace{\sigma_{p-3}\sigma_{p-4} \cdots \sigma_{p-a-1}} \underbrace{\sigma_{p-2}\sigma_{p-3}\cdots \sigma_{p-a-1}}$\\
$\sim_{M}\underbrace{\sigma_{a}\sigma_{a-1}\cdots\sigma_{1}}\underbrace{\sigma_{a+1}\sigma_{a}\cdots \sigma_{2}} \cdots \underbrace{\sigma_{p-3}\sigma_{p-4} \cdots \sigma_{p-a-1}}\underline{\sigma_{p-2}}\underbrace{\sigma_{p-3}\cdots \sigma_{p-a-1}}$\\
$\sim_{M}\underbrace{\sigma_{a}\sigma_{a-1}\cdots\sigma_{1}}\underbrace{\sigma_{a+1}\sigma_{a}\cdots \sigma_{2}} \cdots \underbrace{\sigma_{p-3}\sigma_{p-4} \cdots \sigma_{p-a-1}}\underbrace{\sigma_{p-3}\sigma_{p-4}\cdots \sigma_{p-a-1}}$\\
$\sim_{M}\underbrace{\sigma_{a}\sigma_{a-1}\cdots\sigma_{1}}\underbrace{\sigma_{a+1}\sigma_{a}\cdots \sigma_{2}} \cdots \underbrace{\sigma_{p-3}\sigma_{p-4} \cdots \sigma_{p-a-2}}{(\sigma_{p-3}\sigma_{p-4}\cdots \sigma_{p-a-1})}^2$\\
\begin{center}
$\vdots$
\end{center}
$\sim_{M}\underbrace{\sigma_{a}\sigma_{a-1}\cdots\sigma_{1}}\underbrace{\sigma_{a+1}\sigma_{a}\cdots \sigma_{2}} \cdots\underbrace{\sigma_{p-4}\sigma_{p-5} \cdots \sigma_{p-a-3}}{(\sigma_{p-4}\sigma_{p-5}\cdots \sigma_{p-a-2})}^3 $\\
\begin{center}
$\vdots$
\end{center}
$\sim_{M}\underbrace{\sigma_{a}\sigma_{a-1}\cdots\sigma_{1}}{(\sigma_{a-1}\sigma_{a-2}\cdots \sigma_{1})}^{(p-a-1)} $\\
$\sim_{M}\underline{\sigma_{a}}\underbrace{\sigma_{a-1}\cdots\sigma_{1}}{(\sigma_{a-1}\sigma_{a-2}\cdots \sigma_{1})}^{(p-a-1)} $\\
$\sim_{M}\underbrace{\sigma_{a-1}\sigma_{a-2}\cdots\sigma_{1}}{(\sigma_{a-1}\sigma_{a-2}\cdots \sigma_{1})}^{(p-a-1)} $\\
$={(\sigma_{a-1}\sigma_{a-2}\cdots\sigma_{1})}^{(p-a)} \ = \overline{B}(a,p-a)\sim_M B(a,p-a)$.
\end{pot3}
%
\section{{Minimal Unknotting Crossing Data For Torus Links}}\label{mucdftl}
In this section, we extend the method of finding minimal unknotting crossing data for torus knots to torus links, by small changes in the method.
\begin{thm}\label{thm7}
 Let $K(p, p)$ be a torus link. Then the $U-$crossing data for $B(p, p)$
is a minimal unknotting crossing data for $B(p, p)$ (or $K(p, p)$).
\end{thm}
\begin{pf} It follows from Theorem \ref{thm1}, that $U(B(p,p))$ is a unknotting crossing data for $B(p,p)$. Note that, \[\displaystyle |U(B(p,p))| = \sum_{i=0}^{p-1} i = \frac{(p-1)p}{2},\]
and the unknotting number of $K(p,p)$ is \[ \frac{(p-1)(p-1)+(p,p)-1}{2}=\frac{(p-1)p}{2}.\]
Thus, $U(B(p,p))$ is the minimal unknotting crossing data $B(p,p)$.
\end{pf}
 \begin{figure}
  \begin{center}
   \includegraphics[height=8.2cm,width=13.8cm]{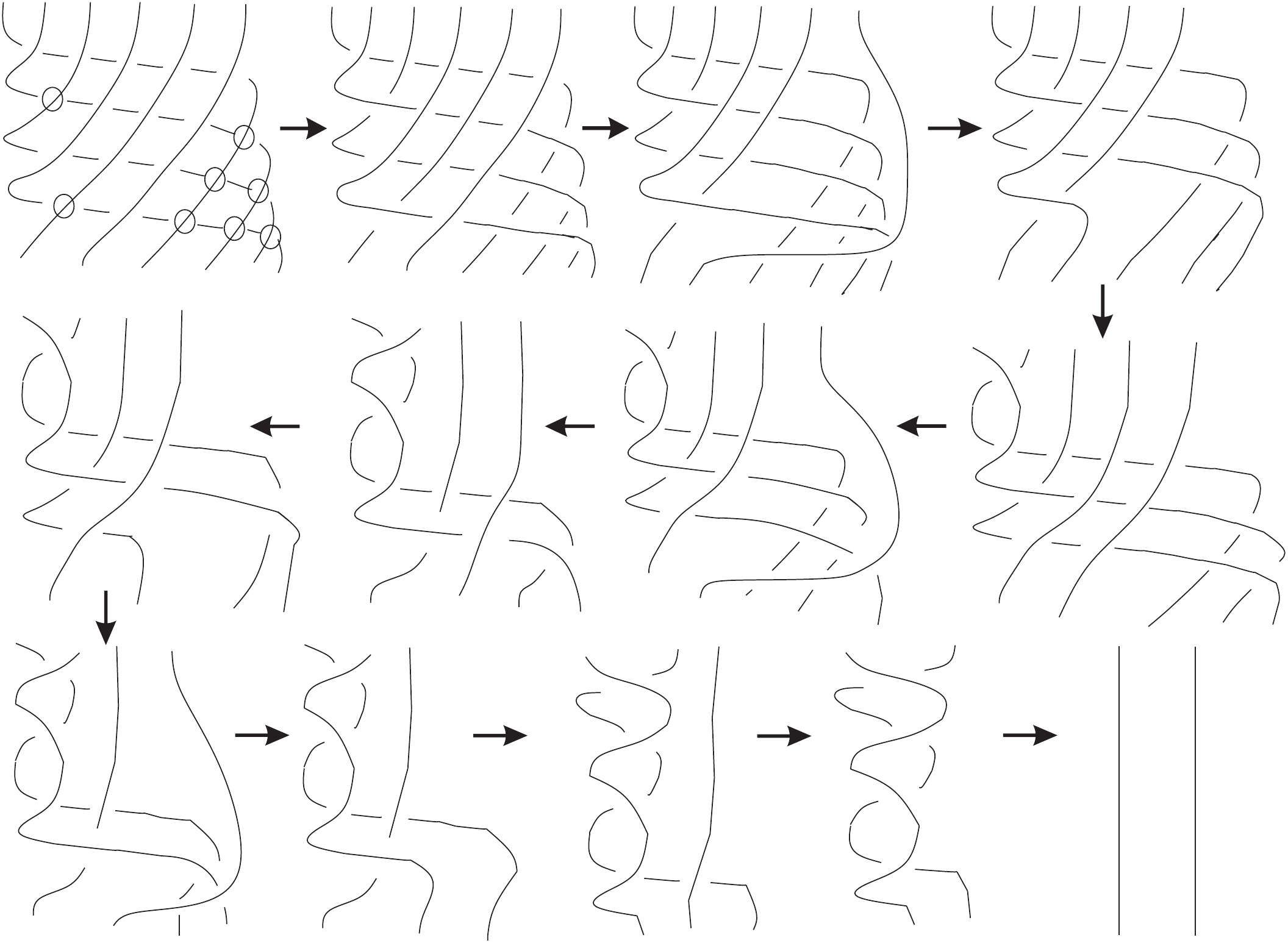}
   \caption{Unknotting procedure for torus link $K(6,4)$}
 \end{center}
\end{figure}
\begin{rmk}\label{rmk6}
 For every $p,q$, with $(p,q)=d \ and\ q \equiv a (mod\ p)$, then $(p,a)=(p-a,a)=d$\ and
\[B(p,q)= {(\sigma_1\sigma_2 \cdots \sigma_{p-2}\sigma_{p-1})}^q={(\sigma_1\sigma_2 \cdots \sigma_{p-2}\sigma_{p-1})}^{mp+a};for\ some\ m\geq 0.\]
Then, after changing the crossings in $B(p,q)$, based on $U(B(p,q))$, we obtain\\
$\alpha ^m \beta $, where
$\alpha=\sigma_1\sigma_2 \cdots \sigma_{p-2}\sigma_{p-1}\sigma_1\sigma_2 \cdots \sigma_{p-2}\sigma_{p-1}^{-1}\cdots \sigma_1^{-1}\sigma_2^{-1} \cdots \sigma_{p-2}^{-1}\sigma_{p-1}^{-1}$
and  $\beta = {\sigma_1\sigma_2 \cdots \sigma_{p-2}\sigma_{p-1}\sigma_1\sigma_2 \cdots \sigma_{p-2}\sigma_{p-1}^{-1}\cdots \sigma_1\sigma_2 \cdots\sigma_{p-(a-1)}^{-1}\sigma_{p-(a-2)}^{-1}\cdots \sigma_{p-2}^{-1}\sigma_{p-1}^{-1}}$.\\
By Theorem \ref{thm4} and Theorem \ref{thm1}, $\alpha $ is a trivial $p$-braid and $\beta$ is Markov equivalent to $B(p-a,a)$.
Thus, we can find a minimal unknotting crossing data for $K(p,q)$, if we can find the unknotting crossing data for
$B(p-a,a)$ such that the number of elements in this unknotting crossing data
is equal to the unknotting number of $K(p, q) \  \smallsetminus\ |U(B(p,q))|$, i.e.,
\[\frac{(p-1)(q-1)+d-1}{2}-\frac{mp(p-1)}{2}-\frac{a(a-1)}{2}
=\frac{(p-a-1)(a-1)+d-1}{2}.\]
Observe that, this is the unknotting number of $K(p-a,a)$.
\end{rmk}

Now we will provide the unknotting process. Consider a torus link of type $K(p,q)$, where $(p,q)=d$ and $q=mp+a$. Let $p_1=p$ and $q_1=q$. Since $(p_1,q_1)=d$, Euclid's algorithm ensures that $\exists\ n\in N$ such that if
we define $q_i,p_i$ as
\[q_{i+1}=q_i-s_i p_i;\ q_{i+1}<p_i\]
\[p_{i+1}=p_i-m_i q_{i+1};\ p_{i+1}<q_{i+1}\]
we get $ q_{n+1}\ or \ p_{n+1}\ =0$.\\
Now, if we select $U$-crossing data for $K(p_i,q_i), K(p_i-q_{i+1},q_{i+1}), K(p_i-2q_{i+1},q_{i+1}),\cdots, K(p_i-(m_i-1)q_{i+1},q_{i+1})$ for $i=1,2,\cdots,n$, then by Theorem \ref{thm4}, Theorem \ref{thm7} and Remark \ref{rmk6}, there are corresponding crossings in $K(p,q)$ which provides minimal unknotting crossing data for $K(p,q)$.
\begin{ex}\label{ex4} By the above said method, minimal unknotting crossing data for $K(6,4)$ is $[6,10,14,15,16,18,19,20]$.
\end{ex}
\textbf{Acknowledgements}\\

Authors thank Professor Akio Kawauchi for his valuable comments and suggestions. Also the first author thanks CSIR, New Delhi and IIT Ropar for providing financial assistance and research facilities.\\


\appendix
\section{Matlab Program For Minimal Unknotting Crossing Data}
\noindent \%To find minimal unknotting crossing data for torus knot $K(p,q)$\\
\%---------------Assume p=p(3) and q=q(3)----------\\
clc\\
clear all\\
p(3) = input(`enter value of p:');\\
q(3) = input(`enter value of q:');\\
\hspace*{.5cm}if $p(3)>q(3)$                                  \hfill \%To find number of digits in q(3)\\
\hspace*{1.0cm}str = num2str(q(3));\\
\hspace*{1.0cm}digits = double(str) - 48;\\ 
\hspace*{1.0cm}nod = sum$(digits>=0)$\\
\hspace*{0.5cm}else                                        \hfill \%To find number of digits in p(3)\\
\hspace*{1.0cm} str = num2str(p(3));\\
\hspace*{1.0cm}digits = double(str) - 48;\\ 
\hspace*{1.0cm}nod = sum$(digits>=0)$                             \hfill \%nod is number of digits in min\{p(3),q(3)\}\\
\hspace*{0.5cm}end\\
\hspace*{0.5cm}a(2) = p(3);\\
\hspace*{0.5cm}a(1) = q(3);\\
\hspace*{0.5cm}V = zeros(nod*5,ceil((p(3)-1)*(q(3)-1)+gcd(p(3),q(3))-1)/2);         \\ \%nod*5 is upper bound for n\\
\hspace*{0.5cm}D = zeros(nod*5,1);\\
\hspace*{0.5cm}count = 3;\\
\hspace*{0.5cm}for i = 3:nod*5\\
\hspace*{1.0cm}count = count+1;\\
\hspace*{1.0cm}if mod(i,2) == 0\\
\hspace*{1.5cm}a(i) = mod(p(i),q(i));\\
\hspace*{1.5cm}m(i) = (p(i)-a(i))/q(i);\\
\hspace*{1.5cm}if $m(i)>=2$\\
\hspace*{2.0cm}B3 = zeros((m(i)-1)*q(i)*(q(i)-1)/2,1);\\
\hspace*{2.0cm}for j = 1:(m(i)-1)\\
\hspace*{2.5cm}for k = 1:(q(i)-1)\\
\hspace*{3.0cm}for g = 1:k\\
B3(((j-1)*(q(i)-1)*q(i)/2)+(k*(k-1)/2)+g,1)=(k*(p(3)-1))+(p(i)-j*q(i)-g);\\
\hspace*{3.0cm}end\\
\hspace*{2.5cm}end\\
\hspace*{2.0cm}end\\
\hspace*{1.5cm}else\\
\hspace*{2.0cm}B3 = [\hspace{.05cm}];\\
\hspace*{1.5cm}end\\
\hspace*{1.5cm}V(i-2,1:size(B3,1)) = B3;           \hfill \%for even i, B(i) is $(i-2)^{th}$ row of V\\
\hspace*{1.5cm}D(i-2,1) = size(B3,1);\\
\hspace*{1.5cm}if mod(p(i),q(i)) $\sim = 1$ $\& \&$  mod(p(i),q(i)) $\sim = $0\\
\hspace*{2.0cm}p(i+1) = a(i);\\
\hspace*{2.0cm}q(i+1) = a(i-1);\\
\hspace*{1.5cm}else break\\
\hspace*{1.5cm}end\\
\hspace*{1.0cm}else                                          \hfill\%if mod(i,2) == 1\\
\hspace*{1.5cm}a(i) = mod(q(i),p(i));\\
\hspace*{1.5cm}m(i) = (q(i)-a(i))/p(i);\\
\hspace*{1.5cm}if m(i)$>=$ 1\\
\hspace*{2.0cm}B1 = zeros((p(i)-1)*p(i)*m(i)/2,1);\\
\hspace*{2.0cm}for j = 1:m(i)\\
\hspace*{2.5cm}for k = 1:(p(i)-1)\\
\hspace*{3.0cm}for g = 1:k\\
B1(((j-1)*(p(i)-1)*p(i)/2)+(k*(k-1)/2)+g,1)=
((j-1)*p(i)+k)*(p(3)-1)+p(i)-g;\\
\hspace*{3.0cm}end \\
\hspace*{2.5cm}end\\
\hspace*{2.0cm}end \\
\hspace*{1.5cm}else \\
\hspace*{2.0cm}B1 = [\hspace{0.05cm}];\\
\hspace*{1.5cm}end\\
\hspace*{1.5cm}if a(i)$>$1 \\
\hspace*{2.0cm}B2 = zeros((a(i)-1)*a(i)/2,1); \\
\hspace*{2.0cm}for k = 1:(a(i)-1)\\
\hspace*{2.5cm}for g = 1:k\\
B2(((k*(k-1))/2)+g,1)=(m(i)*p(i)*(p(3)-1))+(k*(p(3)-1))+(p(i)-g);\\
\hspace*{2.5cm}end\\
\hspace*{2.0cm}end\\
\hspace*{1.5cm}else \\
\hspace*{2.0cm}B2 = [\hspace{0.05cm}];\\
\hspace*{1.5cm}end\\
\hspace*{1.5cm}B = zeros((size(B1,1)+size(B2,1)),1);\\
\hspace*{1.5cm}B = union(B1,B2);                                   \hfill \%B is $U(K(p_i,q_i))$\\
\hspace*{1.5cm}V(i-2,1:size(B,1)) = B\'\ ;                    \hfill\%for odd i, B(i) is $(i-2)^{th}$ row of V\\
\hspace*{1.5cm}D(i-2,1) = size(B,1);\\
\hspace*{1.5cm}if mod(q(i),p(i))$\sim$= 1 \&\& mod(q(i),p(i))$\sim$= p(i)-1 \&\& mod(q(i),p(i))$\sim$=0\\
\hspace*{2.0cm}p(i+1) = a(i-1);\\
\hspace*{2.0cm}q(i+1) = a(i);\\
\hspace*{1.5cm}else break\\
\hspace*{1.5cm}end\\
\hspace*{1.0cm}end\\
\hspace*{0.5cm}end\\
\%-------------To combine all $B(i)^{'s}$----------------\\
\hspace*{0.5cm}n = count-1\\
\hspace*{0.5cm}W = zeros(1,((p(3)-1)*(q(3)-1)+gcd(p(3),q(3))-1)/2);\\
\hspace*{0.5cm}if mod(n,2) == 1\\
\hspace*{1.0cm}W1old = 0;\\
\hspace*{1.0cm}E1old = D(1,1);\\
\hspace*{1.0cm}for zc = 3:2:n\\
\hspace*{1.5cm}W1(zc) = m(zc)*p(zc);\\
\hspace*{1.5cm}W1new = W1(zc)+W1old;\\
\hspace*{1.5cm}W1old = W1new;\\
\hspace*{1.5cm}hh = E1old;\\
\hspace*{1.5cm}E1(zc) = (D((zc-1),1)+D(zc,1));\\
\hspace*{1.5cm}E1new = E1(zc)+E1old;\\
\hspace*{1.5cm}E1old = E1new;\\
\hspace*{1.5cm}WW = (W1old*(p(3)-1))+union(V((zc-1),1:D((zc-1),1)),V((zc),1:D(zc,1)));\\
\hspace*{1.5cm}W(1,1:D(1,1)) = V(1,1:D(1,1));\\
\hspace*{1.5cm}W(1,hh+1:E1new) = WW;\\
\hspace*{1.0cm}end\\
\hspace*{0.5cm}else\\
\hspace*{1.0cm}for zc = 1\\
\hspace*{1.5cm}WW = V(1,1:D(1,1));\\
\hspace*{1.5cm}W(1,1:D(1,1)) = WW;\\
\hspace*{1.0cm}end\\
\hspace*{1.0cm}W1old = 0;\\
\hspace*{1.0cm}E1old = D(1,1);\\
\hspace*{1.0cm}for zc = 3:2:n-1\\
\hspace*{1.5cm}W1(zc) = m(zc)*p(zc);\\
\hspace*{1.5cm}W1new = W1(zc)+W1old;\\
\hspace*{1.5cm}W1old = W1new;\\
\hspace*{1.5cm}jj = E1old;\\
\hspace*{1.5cm}E1(zc) = (D((zc-1),1)+D(zc,1));\\
\hspace*{1.5cm}E1new = E1(zc)+E1old;\\
\hspace*{1.5cm}E1old = E1new;\\
\hspace*{1.5cm}WW = (W1old*(p(3)-1))+ union(V((zc-1),1:D((zc-1),1)),V((zc),1:D(zc,1)));\\
\hspace*{1.5cm}W(1,jj+1:E1new) = WW;\\
\hspace*{1.0cm}end\\
\hspace*{1.0cm}MUKD1 = W\\
\hspace*{1.0cm}MUKD2 = ((p(3)-1)*(q(3)-1)/2)+1-W\\
\hspace*{1.0cm}fprintf(`MUKD1 and MUKD2 are minimal unknotting crossing data')\\
\hspace*{0.5cm}end









\end{document}